# A congruence with the Euler totient function $\varphi$:

Abstract: In this article we give a result obtained of an experimental way for the Euler totient function.

**Theorem 1:** for all $x, y, z, n$ positive integers $x < z$ and $x < y$ we have:

a) $\varphi(x^n + y^n) \equiv 0 \pmod{n}$

b) $\varphi(z^n - x^n) \equiv 0 \pmod{n}$

c) $\varphi(\dfrac{z^n - x^n}{z - x}) \equiv 0 \pmod{n}$

d) $\varphi(\dfrac{x^n + y^n}{x + y}) \equiv 0 \pmod{n}$ in this last case n is odd and $n \geq 5$.

I have verified these results experimentally for n=1 up to 20 and for x, y, z = 1 up to 100.

As particular cases we have that the Cunningham numbers also verify:

$C+(b, n) = b^n + 1$ $\qquad$ $\varphi(C+(b, n)) \equiv 0 \pmod{n}$

$C-(b, n) = b^n - 1$ $\qquad$ $\varphi(C-(b, n)) \equiv 0 \pmod{n}$

and the Mersenne numbers:

$M_n = 2^n - 1$ $\qquad$ $\varphi(M_n) \equiv 0 \pmod{n}$

Also it is verified by the Repunit numbers:

$R(b, n) = \dfrac{b^n - 1}{b - 1}$ $\qquad$ $\varphi(R(b, n)) \equiv 0 \pmod{n}$

I have obtained experimentally also this result that implies the previous theorem:

**Theorem 2:** For all $x, y, n$  $x < y$ positive integers $n \geq 5$ exist $q$ prime number $q \equiv 1 \pmod{n}$  and  $q$ divides $(x^n \pm y^n)/(x \pm y)$.
(n odd in the case +)

Consider: $a = (x^n \pm y^n)/(x+y)$ by the conjecture $a = q^c \cdot b$ with $b$ and c positive integers and $GCD(q,b) = 1$ then we have
$\varphi(a) = \varphi(q^c \cdot b) = \varphi(q^c) \cdot \varphi(b) = q^{c-1}(q-1)\varphi(b)$ and therefore
$\varphi(a) \equiv 0 \pmod{n}$

Proof **[1]**:

One can assume $x$ and y coprime (as m divides n implies phi(m) divides phi (n))
so, they are prime to N= $x^n + y^n$. Now it is clear that a=y/x has order 2n in the group (Z/NZ)*.
(If $a^m = 1$ with m < n, then N divides $y^m - x^m$, which is an integer between 1 and N-1).

Proof **[2]**:

The result is an immediate corollary of the theorem of Birkhoff and Vandiver (Ann Math vol 5, pp. 173-180, 1903).

It states that for integers n > 1, a > b > 0 with a, b coprime then Phi_n(a,b) has a prime factor p not dividing any Phi_m(a,b) for any proper factors m of n EXCEPT if n = 6 and a = 2 and b = 1.

To explain the notation: Phi_n(X,Y) is the homogenized version of the n-th cyclotomic polynomial: Phi_n(X,Y) is the product of X + zY as z runs through the primitive n-th roots of unity.

The prime p whose existence is shown satsifies p = 1 (mod n).

Now $a^n + b^n$ is the product of Phi_k(a,b) for all k dividing 2n but not n, $a^n - b^n$ is the product of Phi_k(a,b) for all k dividing n, for odd n $(a^n + b^n)/(a + b)$ is the product of Phi_k(a,b) for all k>2 dividing 2n but not n, and $(a^n - b^n)/(a-b)$ is the product of Phi_k(a,b) for all k>1 dividing n.

In all cases this number is a multiple of Phi_n(a,b) or Phi_{2n}(a,b).
Assume gcd(a,b) = 1 (the general case reduces to this), unless n = 3 or 6
and a= 2, b = 1, then the number in question has a prime factor
p = 1 (mod n) and phi of it is divisible by p - 1 and so by n.

References:

**[1] : Dominique Bernardi**, Personal communication.
 Theorie des Nombres Institut de Mathematiques - Universite Pierre et Marie Curie
**[2] : Robin Chapman**
Personal communication  www.maths.ex.ac.uk/~rjc/rjc.html

**Author:**

**Sebastian Martin Ruiz**
Avda. de Regla 43 Chipiona 11550 Spain

smruiz(AT)telefonica.net
http://personal.telefonica.terra.es/web/smruiz/